\documentclass[10pt]{amsart}
\usepackage{mathrsfs}
\usepackage{amsfonts}
\usepackage{latexsym,amsmath,amssymb}

 \renewcommand{\epsilon}{\varepsilon}

  \newcommand{\n}{{\mathbf n}}
  \newcommand{\beq}{\begin{equation}}
  \newcommand{\eeq}{\end{equation}}
    \newcommand{\RR}{\mathbb{R}}
    
 \newtheorem{theorem}{Theorem}[section]
 
 \newtheorem{lemma}[theorem]{Lemma}
 \newtheorem{corr}[theorem]{Corollary}
 
 \newtheorem{proposition}[theorem]{Proposition}

 \newtheorem{deff}[theorem]{Definition}
 \newtheorem{rem}[theorem]{Remark}
 \newcommand{\bth}{\begin{theorem}}
 \newcommand{\ble}{\begin{lemma}}
 \newcommand{\bcor}{\begin{corr}}
 \newcommand{\bdeff}{\begin{deff}}
 \newcommand{\bprop}{\begin{proposition}}
 \newcommand{\ele}{\end{lemma}}
 \newcommand{\ecor}{\end{corr}}
 \newcommand{\edeff}{\end{deff}}
 
 \newcommand{\eprop}{\end{proposition}}

 \renewcommand{\Pi}{\varPi}

 \renewcommand{\epsilon}{\varepsilon}

\numberwithin{equation}{section}

\newtheorem{lem}{Lemma}[section]

\title[
The almost global existence for a 3-D  nematic liquid-crystals  equation]
{The almost global existence to classical solution for a 3-D wave equation of nematic liquid-crystals}
\author{Yi Du}
\author{Geng Chen}
\author{Jianli Liu}

\address{Department of Mathematics, South China Normal University,duyidy@gmail.com}
\address{Department of Mathematics,  The Pennsylvania State
University, University Park, PA 16802, USA, chen@math.psu.edu}
\address{Department of Mathematics, Shanghai University, jlliu@shu.edu.cn}
\date{}

\begin{document}
\maketitle

\begin{abstract}
 In this paper we obtain the wave equation modeling the nematic liquid-crystals in three space dimensions and  study the lifespan of classical solution to Cauchy problem. The almost global existence to classical solution for small initial data will be presented.
\end{abstract}


\section{Introduction}
\subsection{Physical background}
Liquid crystals are a state of matter that have properties between those of a conventional
liquid and those of a solid crystal that are optically anisotropic, even when they are at rest.
In this paper, we shall study a variational wave system from liquid crystal. 

At the beginning, we  introduce the physical derivation by the
least action principle for the variational wave equations.
The mean orientation of the long molecules in a nematic liquid crystal is described by a unit vectors 
$\n=\n(t,x)\in \mathbb{S}^2$, the unit sphere, where $(t,x)\in \RR^+\times \RR^3$,
which could be modeled by below Euler-Largrangian equations derived from the least action principle
\cite{AH,[37]},
\beq
\kappa {\mathbf
n}_{tt} + \mu {\mathbf
n}_t + \frac{\delta \int W(\n,\nabla \n)}{\delta \mathbf n} =\lambda {\mathbf
n}, \qquad {\mathbf n}\cdot
{\mathbf n} = 1, \label{1.2}
\eeq
where the well-known Oseen-Franck potential
energy density $W$ is given by
\beq\label{Oseen_Franck}
 \textstyle
W\left({\mathbf n},\nabla{\mathbf
n}\right) = \frac12\alpha(\nabla\cdot{\mathbf n})^2 +
\frac12\beta\left({\mathbf n}\cdot\nabla\times{\mathbf n}\right)^2
+\frac12\gamma\left|{\mathbf n}\times(\nabla\times{\mathbf
n})\right|^2.
\eeq
Here the positive constants $\alpha$,
$\beta$ and $\gamma$ are elastic constants of the liquid crystal,
corresponding to splay, twist, and bend, respectively.
The Lagrangian multiplier $\lambda$ is  determined by the constraint $\mathbf n \cdot \mathbf n=1$.
The inertia and viscous coefficients $\kappa$
and $\mu$ are two non-negative constants, respectively.

There are many studies
on the constrained elliptic system of equations for ${\mathbf n}$,
and on the parabolic flow associated with it, where we refer the reader to see
\cite{[3],[7],[10],[19],[26],[49]}.
In particular, the one-constant elliptic model leads to the equation for harmonic maps taking values in the two-sphere \cite{[7],[10]}.
One parabolic system of equations for ${\mathbf n}$ coupling the compressible Navier-Stokes equations
in one space dimension (1-D) has been considered by \cite{wang}.

However, for the complete hyperbolic system (\ref{1.2}), only the 1-D case has been systematically
studied, where the main difficulty is the possible gradient blowup. For the 1-D Cauchy problem
of (\ref{1.2}) with $\mu=0$ and $\kappa=1$, which describes the model with viscous effects neglected,
the global weak existence and singularity formation
have been extensively studied. In fact, for this extreme case,
one example with smooth initial data and singularity formation (gradient blowup) in finite
time has been provided in
\cite{ghz}, and the global weak existence results have been provides in
\cite{BZ, CZZ12, ZZ03,  ZZ05a, ZZ10, ZZ11}.
In a recent paper \cite{CZ12}, the singularity formation (gradient blowup) and global weak existence for the 1-D Cauchy
problem for the complete system (\ref{1.2}) ($\kappa$ and $\mu$ are non-negative constants) with $\alpha=\gamma$ have been established.
This result shows that for 1-D Cauchy problem we should essentially expect the similar existence and regularity
results for the complete
system (\ref{1.2}) and the extreme case with $\mu=0$ and $\kappa=1$.

Very recently, Li, Witt \& Yin \cite{li-witt-yin} have proved the singularity formation and studied the life-span of the classical solution 
for a 2-D axisymmetric model. 
However, the well-posedness for \eqref{1.2} in 2-D and 3-D is still wide open. 
In this paper, we consider the life-span of the classical solution for the 3-D Cauchy problem of \eqref{1.2} with small initial data.
Note the initial data in the 1-D singularity formation examples in \cite{CZ12, ghz} have large $C^1$-norm.

\subsection{The Variational Wave Equation}
In this paper, we focus on the 3-D Cauchy problem of \eqref{1.2}
with  $\mu=0$ and $\kappa=1$ and planner deformation:
\begin{equation}\label{variation}
\frac{\delta}{\delta \n}\int \bigl\{ (\partial_t \n)^2-
W(\n,\nabla \n)
\bigr\}\,dx\, dt=0, \quad |\n|=1,
\end{equation}
where the planer deformation means that
 \beq\label{functions of n}
\n=(n_1,n_2,n_3)=(\cos u(t,x), \sin u(t,x),0)
 \eeq
with the angel $u=u(t,x)$,  $(t,x)\in \RR^+ \times \RR^3$.

Now, the variation \eqref{variation} can be reformulated as
 \begin{equation}\label{variation1}
\frac{\delta}{\delta \n}\int\bigl\{ \frac12(\n_t)^2-
W(\n,\nabla \n)\bigr\}\,dx\, dt =0,
\end{equation}
which gives rise to the Euler-Lagrange equation
 \begin{equation}\label{Euler-Lagarange}
\partial_{tt}n_i+\partial_{n_i}W(\n,\nabla \n)-\sum_{j=1}^3\partial_{x_j}\bigl(\partial_{\partial_j n_i}W(\n,\nabla \n)\bigr)=0.
\end{equation}
Multiplying three equations in \eqref{Euler-Lagarange} by $-\sin u$ and $\cos u$,  then adding them up,
we can get a variational wave equation, by (\ref{functions of n}) and (\ref{Oseen_Franck}),
\begin{equation}\label{original wave equation}
u_{tt}-c_1^2(u)\,\partial_{x_1}^2u -c_2^2(u)\,\partial_{x_2}^2u -\beta\, \partial_{x_3}^2 u + 
(\alpha-\gamma)\,\sin2u\, \partial_{x_1}\partial_{x_2} u=F(u,\partial u),
\end{equation}
where
\[
c_1(u)^2=\alpha\sin^2u +\gamma\cos^2u,\quad c_2^2(u)=\alpha\cos^2u+\gamma\sin^2u,
\]
and
 \begin{equation}\label{F}
	F(u,\partial u)=\frac12(\alpha-\gamma)\sin2u\,[(\partial_1u)^2-(\partial_2u)^2]-(\alpha-\gamma)	
	\cos2u\,\partial_1u\,\partial_2u.
 \end{equation}
Instead of the direct calculation from (\ref{Euler-Lagarange}), another way to derive (\ref{original wave equation}) is 
from the variational principle
\[
\frac{\delta}{\delta u}\int \bigl\{ (\partial_t \n)^2-
W(\n,\nabla \n)\bigr\}\,dx\, dt
=0, \quad |\n|=1,
\]
where in the derivation we need to use (\ref{functions of n}) and 
\begin{equation}\label{function W}
W=\frac{\alpha}{2}(\cos u\,\partial_{x_2} u-\sin u\,\partial_{x_1}u)^2+
\frac{\beta}{2}(\partial_{x_3}u)^2+\frac{\gamma}{2}(\cos u\,\partial_{x_1}u+\sin u\,\partial_{x_2} u)^2.
\end{equation}
Since $\alpha, \ \beta,$ and $\gamma$ are all positive constants,
we can rewrite \eqref{original wave equation} as the following  equation by scalling
\begin{equation}\label{0}
u_{tt}-\sum\limits_{i,j=1}^3a_{ij}(u(t,x))\partial_{i}\partial_j u=F(u,\partial u(t,x)), \qquad x\in \RR^3,\quad
t\in \RR^+,
\end{equation}
where $\partial_i=\frac\partial{\partial x_i}, \ (i=1,2,3)$ and
\begin{equation}
A=(a_{ij})_{1\leq i,j\leq 3}=\left(
  \begin{array}{ccc}
    c_1^2(u)& \frac12(\alpha-\gamma)\sin 2u  & 0 \\
    \frac12(\alpha-\gamma)\sin 2u & c_2^2(u) & 0 \\
    0 & 0 & 1 \\
  \end{array}
\right).
 \end{equation}
It is obvious that there exist  constants $M_0,M_1$, such that
\[
0<M_0\leq c_1^2(u),\,c_2^2(u)\leq M_1.
\]
Furthermore,  $A$ is a positive definite symmetric matrix.

In our paper, we shall study the lifespan for the variational wave equation \eqref{original wave equation} in three dimensional space with small and smooth initial data
\begin{equation}\label{initial data}
u(0,x)=\epsilon_0f(x),\ u_t(0,x)=\epsilon_0 g(x),
\end{equation}
where $\epsilon_0$ is a constant small enough, and $(f,g)\in C^\infty (\RR^3)$ with compact support.
    We will obtain the following almost global existence result:
\begin{theorem}\label{main}
The equations \eqref{original wave equation} and  \eqref{initial data} or equivalently  \eqref{0} and  \eqref{initial data}
admit a unique solution in $C^\infty([0,T^\star), \RR^3)$ with
\begin{equation}
T^\star\geq C\exp {\frac{1}{\epsilon_0}},
\end{equation}
where C is a constant independent of $\epsilon_0$.
\end{theorem}
We have already known from \cite{CZ12, ghz}
that there are 1-D gradient blowup examples when the initial data have large $C^1$
norm, even when the amplitude oscillation of the initial data is small. We note that the 1-D examples
are also special examples in 3-D with initial data restricted in one space dimension.
More precisely, for the blowup examples in  \cite{CZ12, ghz}, in which they also consider the
planar transformation $\n=(\cos u, \sin u,0)$, 
the initial data $u_0=u(0,x)$ satisfy $\|u_0\|_{L^\infty}=O(\epsilon)$, $\|u_0\|_{C^1}=O(1)$
and the blow-up time is at $O(1)$, where $\epsilon>0$ is an arbitrarily small number.
Clearly, these singularity formation results of  \cite{CZ12, ghz}  do not conflict with our result.
On the contrary, combining these two pieces of information from \cite{CZ12, ghz}  and this paper,
we have a fairly complete picture of the smooth classical solutions  for the 3-D liquid crystal equation  (\ref{0}).

We remark for the axisymmetric case.  In \cite{li-witt-yin}, the authors have proved the singularity formation and 
 the lifespan for the axisymmetric case in 2-D with small data.   Due to the significant difference between 2-D and 3-D for the wave equation, 
for the 3-D axisymmetric case,
whether the solution of the equations \eqref{original wave equation} with small initial data form singularity or not is still open,  
although  the nonlinear structure  of our model is similar to the one in
Lindblad \cite{lindblad1992}. We leave this case to be concerned in the future paper.

To prove Theorem \ref{main}, we use the generalized energy method which is popular  in analyzing wave equations.
There are large amount of literatures focusing on this topic. See \cite{hormander, klainerman85, klainerman87,  li-yu, li-zhou1992, lindblad1990, lindblad1992}
and \cite{li-chen} for more details.  
In Lei-Lin-Zhou \cite{lei-lin-zhou}, the authors have presented a global existence for Faddeev model 
in both 2-D and 3-D cases, where the equations have a null condition structure.
However, the equation \eqref{0} does not have the null condition structure, so we can not get the
same results.


This paper is organized as following: in the next section, we recall some well-known results for the wave equation, and in Section 3, we shall prove our main results.

 \section{Preliminary}
Without loss of generality, in the following we assume the positive constants $\alpha$ and $\gamma$ satisfy
\[
\alpha\leq \gamma.
\]
 The proof for the case $\alpha> \gamma$ is very similar.
Then, by a schedule of re-scaling, we can rewrite the \eqref{0} in the following wave equation
 \begin{equation}\label{2.1}
\begin{cases}
\Box u-\sum\limits_{i,j=1}^3\bar{a}_{ij}(u(t,x))\partial_{i}\partial_j u=\bar{F}(u,\partial u(x,t)), \qquad x\in \RR^3,\
t>0,\\
t=0,u=\epsilon f, u_t=\epsilon g,
\end{cases}
\end{equation}
here and hereafter 
\begin{equation}\Box =\partial_{t}^2-\sum\limits_{i=1}^3\partial_{x_i}^2,
\end{equation} 
and
\begin{equation}
\bar{A}=(\bar{a}_{ij})_{1\leq i,j\leq 3}=\left(
  \begin{array}{ccc}
    c_0\cos^2u& \frac12c_0\sin2u & 0 \\
    \frac12c_0\sin2u & c_0\sin^2u & 0 \\
    0 & 0 &0 \\
  \end{array}
\right)
 \end{equation}
 is  positive semi-definite with $c_0=\frac{\gamma-\alpha}{\alpha}\geq0$.

In the beginning, we introduce the following
notations. Write
\beq\left\{
\begin{array}{rcl}
x_0&=&t,\\
\partial_0&=&-\partial/\partial_t,\\
\partial_i&=&\partial/\partial_{x_i},\\
\partial&=&(\partial_0,\partial_1,\partial_2,\partial_3)=(-\partial_t,\partial_1,\partial_2,\partial_3).
\end{array}
\right.
\eeq
Define the first order differential
operators as
 \begin{equation}
 \Omega =(\Omega_{ab})_{0\leq a,b \leq3},
\end{equation}
with
 \begin{equation}
 \Omega_{ab}=x_a\partial_b-x_b\partial_a,\ (0\leq a,b \leq3),
\end{equation}
and the scaling operator as
\begin{equation}
L=\sum\limits_{i=0}^3x_i\partial_i.
\end{equation}
Denote the vector field
\begin{equation}
\Gamma=(\partial,\Omega,L).
\end{equation}
It is well known that the above operators $\Gamma$ have possessed a good commutation with the wave operator $\Box=\partial_{t}^2-\Delta $ (See also \cite{klainerman85, li-chen}):
\begin{lem}
 For multi-indices  $\varsigma,\
 \xi$,
we have
\begin{equation}
[\Box ,\Gamma^\xi]= \sum_{|\varsigma|\leq
|\xi|-1}A_{\xi\varsigma}\Gamma^\varsigma\Box,
\end{equation}
where $[\ ,\ ]$ stands for poisson bracket and $A_{\xi\varsigma}$ are
constants.
\end{lem}


For any integer $N\geq0,$  denote
\begin{equation}
\| v(t,\cdot)\|_{\Gamma,N,p} =\sum_{|\varsigma|\leq
N}\|\Gamma^\varsigma v(t,\cdot)\| _{L^{p}}.
\end{equation}
Then, we have the following decay estimate

\begin{lem}\label{lem_2.2}
Suppose that $h = h(t,x)$ with $(t,x)\in \RR^+\times \RR^n$ is a function with compact support
in the variable $x$ for any fixed $t \geq 0$. For any integer $N \geq 0$,we have
\begin{equation}
\| h(t,\cdot)\|_{\Gamma,N,\infty} \leq C(1+t)^{-\frac{2}{p}}(1+|t-|x||)^{-\frac1p}
\|h(t,\cdot)\| _{\Gamma, N+[\frac{n}{p}]+1,p},
\end{equation}
where $p\geq1$, C is a constant.
\end{lem}
\begin{rem}
The proof of Lemma 2.2 can be found in \cite{li-zhou1995}.
\end{rem}

In order to prove Theorem \ref{main}, we also need the following estimate for a linear wave equation (See \cite{li-zhou1992,lindblad1990}):
\begin{lem} \label{lemma_2.3}
Let $h(t,x)$ and $w(t,x)$ be the functions with support $\{x||x|\leq t+\rho\}$ in the variable $x$
for any fixed $t\geq 0$,  if for all norms on the right hand side of below inequality are bounded, then holds
\begin{equation}
  \|h\partial w\|_{L^2}\leq C_\rho \|\nabla h\|_{L^2}\sum\limits_{|\alpha|\leq 1}\|\Gamma^\alpha w\|_{L^\infty}.
\end{equation}
\end{lem}
\section{the almost global existence}
In this section, we will prove Theorem \ref{main}. Denote
\begin{equation}\label{3.1}
 E_{s,T}= \{v=v(x,t)|   \sup\limits_{0\leq t\leq T}\|\partial v\|_{\Gamma, s,2}\leq \epsilon,  \mbox{ with} \  v(0,x)=\epsilon_0 f, v_t(0,x)=\epsilon_0 g\}.
\end{equation}
where $s, \epsilon_0$ and $\epsilon$ are  given positive constants with  $s>9$ and $\epsilon_0$, $ \epsilon$ small enough.
We shall prove our main result by finding a fixed point in $E_{s, T}$. First, we give the linearized equation of \eqref{0} as the following case:
\begin{equation}\label{3.2}
\begin{cases}
\Box u-\sum\limits_{i,j=1}^3\bar{a}_{ij}(v(t,x))\partial_{i}\partial_j u=\bar{F}(v,\partial v(x,t)), \qquad x\in \RR^3,\
t>0,\\
t=0,u=\epsilon_0 f, u_t=\epsilon_0 g,
\end{cases}
\end{equation}
where  $v\in E_{s,T}$ is any given function. Therefore, the  equation \eqref{3.2}
is a linear wave equation, which admits a unique global solution. We define a
map 
\[Mv=u.\] 
It is sufficient  to prove the map $M$ is a contract map. To do this, we apply the operator $\Gamma^k$ on both side of  \eqref{3.2} (here and  hereafter  k is a multi-indices with  $|k|=s$). By  Lemma 2.1, we have
 \begin{align}\label{3.3}
 &\Box\Gamma^k u-\sum\limits_{i,j=1}^3\bar{a}_{ij}(v(t,x))\partial_{i}\partial_j\Gamma^k u\\\nonumber
  &=
  \sum\limits_{|\alpha|\leq |k|}C_{k\alpha}\Gamma^\alpha\bar{F}(v,\partial v(x,t)) +\sum\limits_{i,j=1}^3[\bar{a}_{ij}(v(t,x)\partial_{i}\partial_j,\Gamma^k] u,
 \end{align}
where $[\ ,\ ]$ is the poisson bracket.
Multiplying $\partial_t\Gamma^k  u$ and taking integral on $\RR^3$ for \eqref{3.3}, we get
 \begin{align}\label{3.4}
 &\frac12\frac{d}{dt}\bigg(\|\Gamma^k\partial u\|_{L^2(R^3)}^2+ \sum\limits_{i,j=1}^3\int_{R^3}\bar{a}_{ij}(v(t,x))\partial_{i}\Gamma^k u_i\partial_j\Gamma^k u_idx\bigg) \\\nonumber
  =&
  \sum\limits_{|\alpha|\leq |k|}C_{k\alpha}\int_{R^3}\Gamma^\alpha\bar{F}(v,\partial v(x,t))\partial_t \Gamma^\alpha u dx &
  \\\nonumber
  &+\sum\limits_{i,j=1}^3\int_{R^3}[\bar{a}_{ij}(v(t,x)\partial_{i}\partial_j,\Gamma^k] u\partial_t\Gamma^k  u dx
  \\\nonumber
  &+\frac12\int_{R^3}\partial_t \bar{a}_{ij}(v(t,x))\partial_i\Gamma^k u_i  \partial_j\Gamma^k u_i dx
  \\\nonumber
  &-\frac12\int_{R^3}\partial_j \bar{a}_{ij}(v(t,x))\partial_i\Gamma^k u_i  \partial_t\Gamma^k u_i dx,
 \end{align}
 where $C_{k\alpha}$ are constants. Then, we get
 \begin{align}\label{3.5}
 & \|\Gamma^k\partial u\|_{L^2(R^3)}^2+ \sum\limits_{i,j=1}^3\int_{R^3}\bar{a}_{ij}(v(t,x))\partial_{i}\Gamma^k u_i\partial_j\Gamma^k u_idx  \\\nonumber
 \leq & C\epsilon_0 +C\bigg(
  \sum\limits_{|\alpha|\leq |k|}
  \int_0^t\int_{R^3}\Gamma^\alpha\bar{F}(v,\partial v(x,t))\partial_t \Gamma^\alpha u dxdt &
  \\\nonumber
  &+\sum\limits_{i,j=1}^3\int_0^t\int_{R^3}[\bar{a}_{ij}(v(t,x),\Gamma^k] \partial_{i}\partial_j u\partial_t \Gamma^k u_i dxdt
  \\\nonumber
  &+\frac12\int_{R^3}\int_0^t\partial_t \bar{a}_{ij}(v(t,x))\partial_i\Gamma^k u_i \partial_j \Gamma^k u_i dxdt
  \\\nonumber
  &-\frac12\int_{R^3}\int_0^t\partial_j \bar{a}_{ij}(v(t,x))\partial_i\Gamma^k u_i\partial_t \Gamma^k  u_i dxdt
  \bigg)
  \\\nonumber
  =&C \epsilon_0 + I + II + III + IV.
 \end{align}
By the chain rule and lemma 2.2, we get
\begin{align}\label{3.6}
 III+IV\leq  & C\|\partial u\|_{\Gamma, s,2}^2\sum\limits_{i,j=1}^3\int_0^t \|\partial \bar{a}_{ij}(v(t,x))\|_{L^\infty}dt \\\nonumber
  \leq  & C \|\partial u\|_{\Gamma, s,2}^2\int_0^t \|\partial v\|_{L^\infty}dt
  \leq   C\epsilon\|\partial u\|_{\Gamma, s,2}^2\int_0^t<\tau>^{-1}d\tau,
\end{align}
here and hereafter, we use the notation
\begin{equation}\nonumber
<x>=\sqrt{1+x^2}.
\end{equation}
Recalling $s>9$ and $|k|=s$,  we have,
\begin{align}\label{3.7}
 II\leq  &C \sum\limits_{i,j=1}^3\int_0^t\int_{R^3} \sum\limits_{|\beta|\leq |k|-1 }\Gamma^\beta \big(\bar{a}_{ij}(v(t,x))\partial_{ij}u\big)\partial \Gamma^ku dxdt
 \\\nonumber
 \leq  &C \sum\limits_{i,j=1}^3\int_0^t\int_{R^3} \sum\limits_{|\alpha|\leq [\frac{|k|-1}{2}] } \sum\limits_{|\beta|\leq |k|-1 }\big(|\Gamma^\alpha \bar{a}_{ij}(v(t,x))\partial_{ij} \Gamma^\beta u|\\\nonumber
 & +| \partial_{ij}\Gamma^\alpha u\Gamma^\beta \bar{a}_{ij}(v(t,x)) |\big) | \partial\Gamma^k u| dxdt\\\nonumber
  \leq & C\int_0^t \|v\|_{\Gamma,[\frac{|k|}{2}],\infty}\|\partial u\|_{\Gamma,|k|,2}^2 +\|u\|_{\Gamma,[\frac{|k|}{2}],\infty}\|\partial u\|_{\Gamma,|k|,2}  \|\partial v\|_{\Gamma,|k|,2} d\tau\\\nonumber
 \leq & C\epsilon\big(|\partial u\|_{\Gamma,s,2}+1\big)\|\partial u\|_{\Gamma,s,2}\int_0^t<\tau>^{-1}d\tau,
\end{align}
as well as
 \begin{align} \label{3.8}
   I\leq & C\int_0^t\big[\sum\limits_{|\alpha|\leq [\frac{|k|}{2}]}(\|\Gamma^\alpha \sin 2v\|_{L^\infty}+\|\Gamma^\alpha\partial v\|_{L^\infty})
  \\ \nonumber
  &\cdot \sum\limits_{|\beta|\leq |k|}(\|\Gamma^\beta (\sin 2v\partial v)\|_{L^2}+\|\Gamma^\beta(\partial v)^2\|_{L^2})\big]dt\cdot \|\partial u\|_{\Gamma, |k|, 2}
  \\ \nonumber
\leq& C\int_0^t
  \bigg(\|v\|_{\Gamma, [\frac{k}2]+1,\infty}+ O\big(\|v\|_{\Gamma, [\frac{k}2]+1,\infty}^3\big)\bigg)^2
  \|\partial v\|_{\Gamma,|k|,2}dt \cdot \|\partial u\|_{\Gamma, |k|,2}
  \\ \nonumber
  \leq&  C\epsilon^2 \|\partial u\|_{\Gamma, s, 2}.
  \end{align}
where we used the Taylor's formula.
 Recalling that $(\bar{a}_{ij})$ is positive semi-definite, then by (\ref{3.4})-(\ref{3.8}), we get
\begin{align}\label{3.9}
  \|\partial u\|_{\Gamma, s,2}^2\leq C\epsilon_0 +C'\epsilon  \|\partial u\|_{\Gamma, s,2}^2 \ln(1+t),
     \end{align}
     where $C$ and $C'$ are constants. Then from \eqref{3.9}, let $T\leq T_0\leq e^{\frac1{2C'\epsilon}}$ and we take the initial data small enough such that 
\[C\epsilon_0<\frac14\epsilon\ll1,\]  
then 
\[Mv=u\in E_{s,T}.\] 
One still needs to prove $M$ is a contraction map.
To do this, let $v_1, v_2\in E_{s,T}$, and then  
\[Mv_i=u_i\in E_{s,T}, \quad i=1,2.\]  
Subsequently, we need to estimate
$\|\partial u_1-\partial u_2\|_{\Gamma,s,2}$. For simplicity, we write $U=u_1-u_2$, and correspondingly, $V=v_1-v_2$,  then from equation \eqref{2.1}, we have
\begin{align}\label{3.10}
\Box U =& \sum\limits_{i,j=1}^3\bigg(\bar{a}_{ij}(v_1(t,x))\partial_{ij}u_1-\bar{a}_{ij}(v_2(t,x))\partial_{ij}u_2\bigg)\\\nonumber
&+F(v_1,\partial v_1)-F(v_2,\partial v_2)  \\\nonumber
  =&
\sum\limits_{i,j=1}^3\bigg(\bar{a}_{ij}(v_1(t,x))\partial_{ij}U\bigg)+\sum\limits_{i,j=1}^3\bigg([\bar{a}_{ij}(v_1)-\bar{a}_{ij}(v_2)]
\partial_{ij}u_2\bigg)
\\\nonumber
&+F(v_1,\partial v_1)-F(v_2,\partial v_2).
\end{align}
Then by the standard energy estimates, 
we get
\begin{align}\label{3.11}
&\|\partial U\|_{\Gamma,s,2}+ \sum\limits_{i,j=1}^3\int_{R^3}\bar{a}_{ij}(v_1(t,x))\partial_{i}\Gamma^\alpha U\partial_j\Gamma^\alpha Udx\\\nonumber
\leq & \sum\limits_{|\alpha|\leq |k|}\sum\limits_{i,j=1}^3 \int_0^t\int_{R^3} |[ \bar{a}_{ij}(v_1)\partial_{ij}, \Gamma^\alpha ]U \partial_t\Gamma^k U  |dxdt  \\\nonumber
&+
\sum\limits_{i,j=1}^3 \int_0^t\int_{R^3} |\partial \bar{a}_{ij}(v_1)||\partial\Gamma^k U|^2dxdt  \\\nonumber
 & +\sum\limits_{i,j=1}^3\int_0^t\int_{R^3}|\Gamma^k\big([\bar{a}_{ij}(v_1)-\bar{a}_{ij}(v_2)]\partial_{ij}u_2\big)\partial_t \Gamma^k U| dxdt
\\\nonumber
&+\int_0^t\int_{R^3}|\Gamma^k \big(F(v_1,\partial v_1)-F(v_2,\partial v_2)\big)\partial_t \Gamma^kU|dxdt  \\\nonumber
=&
I+II+III+IV,
\end{align}
where $|k|=s$. Recalling that $v_1, v_2\in E_{s,T}$ and using Lemma \ref{lem_2.2}, we can do similar estimates as  (\ref{3.4})-(\ref{3.8})  
\begin{align}\label{3.12}
I+II\leq & C\|\partial U\|_{\Gamma,s,2}^2\int_0^t \sum\limits_{i,j=1}^3\|\partial \bar{a}_{ij}(v_1(t,x))\|_{L^\infty}dt\\\nonumber
&+C\sum\limits_{i,j=1}^3\int_0^t\int_{R^3} \sum\limits_{|\alpha|\leq [\frac{|k|-1}{2}] } \sum\limits_{|\beta|\leq |k|-1 }\big(|\Gamma^\alpha \bar{a}_{ij}(v_1(t,x)) \partial_{ij}\Gamma^\beta U|\\\nonumber
 & \qquad +|\partial_{ij}\Gamma^\alpha  U\Gamma^\beta \bar{a}_{ij}(v_1(t,x)) |\big) |\partial\Gamma^k  U| dxdt\\\nonumber
\leq & C''\|\partial U\|_{\Gamma,s,2}^2 \ln(1+t).
\end{align}
where we used Taylor's formula to deal with $\bar{a}_{ij}(v_1)$, and hereafter,  $C''$ is a uniform constant.   Doing the expansion
\begin{equation}\label{3.13}
\big(\bar{a}_{ij}(v_1)-\bar{a}_{ij}(v_2)\big)=\left(
                                      \begin{array}{ccc}
                                        2\sin(v_1+v_2)\sin V & 2\cos(v_1+v_2)\sin  V  & 0\\
                                         2\cos(v_1+v_2)\sin V &  \sin(v_1+v_2)\sin V & 0 \\
                                        0 & 0 & 0 \\
                                      \end{array}
                                    \right),
\end{equation}
we can deal with $III$ as following
\begin{align}\label{3.14}
 III \leq  &C\sum\limits_{i,j=1}^3\int_0^t\int_{R^3} \sum\limits_{|\alpha|\leq [\frac{|k|}{2}] } \sum\limits_{|\beta|\leq |k|}\bigg(|\Gamma^\alpha \sin(V(t,x)) \partial_{ij}\Gamma^\beta u_2|\\\nonumber
 & +|\partial_{ij}\Gamma^\alpha  u_2\Gamma^\beta \sin(V(t,x)) |\bigg) |\partial \Gamma^k U| dxdt\\\nonumber
  \leq & C\int_0^t \|\partial V\|_{\Gamma,[\frac{|k|}{2}],\infty}\|\partial u_2\|_{\Gamma,|k|,2} +\|u_2\|_{\Gamma,[\frac{|k|}{2}],\infty}\|\partial V\|_{\Gamma,|k|,2} d\tau \|\partial V\|_{\Gamma,|k|,2}\\\nonumber
 \leq & C''\epsilon \big(\|\partial V\|_{\Gamma,s,2}^2+\|\partial U\|_{\Gamma,s,2}^2\big)\ln(1+t),
\end{align}
where we use Lemma \ref{lemma_2.3} to get the second inequality. Similarly,
we can deal with $IV$ as
\begin{align}\label{3.15}
 IV \leq  &C \sum\limits_{i,j=1}^3\int_0^t\int_{R^3} \sum\limits_{|\alpha|\leq [\frac{|k|}{2}] } \sum\limits_{|\beta|\leq |k|}\bigg(|\Gamma^\alpha \big(\cos2v_1(t,x)\partial v_1 \big)\partial\Gamma^\beta V|\\\nonumber
 & \qquad +|\partial\Gamma^\alpha  V\Gamma^\beta \big(\cos2v_1\partial v_1\big)|+|\Gamma^\alpha \big(\partial v_2\cos2(v_1+v_2)\big)\Gamma^\beta\big(\partial v_2\cos V\big)|\\\nonumber
 &+|\Gamma^\alpha \big(\partial v_2 \cos V\big)\Gamma^\beta\big(\partial v_2 \cos(v_1+v_2)\big)|\bigg)|\partial\Gamma^k  U| dxdt\\\nonumber
   \leq & C\epsilon \ln(1+t)\bigg(1+\int_0^t
    \big(\| V\|_{\Gamma,[\frac{|k|}{2}+1],\infty}^2+\|V\|_{\Gamma,[\frac{|k|}{2}+1],\infty}^4\big)d\tau \\\nonumber
  & +C\int_0^t\epsilon <\tau>^{-1}\big(1+\epsilon<\tau>^{-1}+<\tau>^{-1}\|\partial V\|_{\Gamma,|k|,2}^2\big)d\tau \bigg)\|\partial V\|_{\Gamma,|k|,2}^2\\\nonumber
\leq  & C''\epsilon \ln(1+t) \|\partial V\|_{\Gamma,s,2}^2+C''<t>^{-1} \|\partial V\|_{\Gamma,s,2}^4\,.
\end{align}
Then from (\ref{3.11})- (\ref{3.15}), we have
\begin{align}\label{3.16}
&\|\partial U\|_{\Gamma,s,2}^2+ \sum\limits_{i,j=1}^3\int_{R^3}\bar{a}_{ij}(v_1(t,x))\partial_{i}\Gamma^\alpha U\partial_j\Gamma^\alpha Udx\\\nonumber
\leq &C''\epsilon \ln(1+t) \|\partial V\|_{\Gamma,s,2}^2+C''\epsilon \ln(1+t) \|\partial U\|_{\Gamma,s,2}^2+C''<t>^{-1} \|\partial V\|_{\Gamma,s,2}^4\,.
\end{align}
Let  $T\leq T_1\leq e^{\frac1{4C'' \epsilon}}$, then there holds
\begin{equation}\label{3.17}
\|\partial U\|_{\Gamma,s,2}^2+ \sum\limits_{i,j=1}^3\int_{R^3}\bar{a}_{ij}(v_1(t,x))\partial_{i}\Gamma^\alpha U\partial_j\Gamma^\alpha Udx
\leq \frac34 \|\partial V\|_{\Gamma,s,2}^2\,.
\end{equation}
Combining (\ref{3.9}) and (\ref{3.17}),  take $T\leq \min\{T_0, T_1\}$, we get the desired result.
\bigskip

{\bf ACKNOWLEDGMENTS}\quad  The first author is partially supported by the
NSFC (No. 11001088) and the pearl river new star (No. 2012001). The third author is partially supported by the NSFC-Tianyuan Special Foundation
(No.11126058), Excellent Young Teachers Program of Shanghai
         and the Shanghai Leading Academic Discipline Project
(No. J50101). This work was done when Yi Du and Jianli liu were visiting the
Department of Mathematics of Penn State University during 2012.
They would like to thank professor Chun Liu
and the institute for their the hospitality.

%

\end{document}